%% file: main.tex
\documentclass{article}
\usepackage{makeidx} 
\usepackage[utf8]{inputenc}
\usepackage{hyperref}
\usepackage{colortbl}
\usepackage[margin=1.25in]{geometry}
\usepackage[normalem]{ulem}
\usepackage{lscape}
\usepackage{wrapfig}
\usepackage{booktabs}
\usepackage{bbm}
\usepackage{xcolor}
\usepackage{graphicx}
\usepackage{amssymb}
\usepackage{amsbsy}
\usepackage{amsmath}
\usepackage{amsthm}
\usepackage{amsfonts}
\usepackage{enumitem}
\usepackage{subcaption}
\usepackage{caption}
\usepackage{listings}
\lstdefinelanguage{Julia}{
  keywords={function,end,if,else,elseif,for,while,return,true,false,
             begin,let,try,catch,finally,do,in,import,using,export,
             abstract,type,immutable,mutable,struct,module,baremodule,
             macro,quote,local,global,const,nothing},
  sensitive=true,
  comment=[l]{\#},
  morecomment=[s]{\#=}{=\#},
  morestring=[b]",
  morestring=[b]',
}
\usepackage{authblk}
\usepackage[]{natbib}
\usepackage{enumitem}
\usepackage{algorithm}
\usepackage{tabularx}
\usepackage{longtable}
\usepackage{algpseudocode}
\usepackage{accents}
\usepackage{multicol}
\usepackage{pgfplots}
\pgfplotsset{compat=1.18}
\usetikzlibrary{patterns, arrows.meta}
\usepgfplotslibrary{fillbetween}

\newcommand{\ignore}[1]{}

\hypersetup{
    colorlinks=true,
    linkcolor=blue,
    filecolor=magenta,      
    urlcolor=cyan,
}
\DeclareMathOperator{\relint}{relint}

\DeclareMathOperator{\dom}{dom}
\DeclareMathOperator{\epi}{epi}
\DeclareMathOperator{\conv}{conv}

\newcommand{\ti}[1]{\tilde{#1}}

\newcommand{\un}[1]{\underline{#1}}

\newcommand{\mbf}[1]{\mathbf{#1}}

\title{SDDmiP.jl: A Software Package with a Provably Convergent Benders Algorithm for Multi-Stage Stochastic Mixed-Integer Programming}
\date{\today}
\author{Akul Bansal}
\author{Simge K\"u\c{c}\"ukyavuz}
\affil{Industrial Engineering and Management Sciences \\ Northwestern University \\
{\texttt{akulbansal5@gmail.com, simge@northwestern.edu}}}

\begin{document}
\maketitle

\begin{abstract}
    We present an open-source software package that implements a provably convergent Benders-type decomposition algorithm for multistage stochastic integer programs. In addition to standard cut families, such as Benders, strengthened Benders, and Lagrangian cuts, the algorithm incorporates rectified linear unit (ReLU) cuts, which provide convergence guarantees for general mixed-integer state variables. However, the dual problems used to generate these cuts often admit multiple optimal solutions. Although each solution yields a valid cut that separates the incumbent, the resulting cuts can differ in how well they approximate the subproblem cost. To strengthen these cuts, our package implements and evaluates two cut-selection strategies based on normalization and regularization of the dual problem. We also incorporate an alternating-cut criterion that uses cheaper Benders cuts when they are effective and invokes more expensive tight cuts only when necessary. Computational experiments on four classes of multistage stochastic integer programs benchmark these methods and provide insights on how problem structure affects their practical performance.
\end{abstract}


\section{Introduction}

A multistage stochastic integer program (MSIP) is a framework for sequential decision-making under uncertainty, in which decisions at each stage must be made before the realization of future uncertainty, while satisfying integrality constraints on certain decision variables. The uncertainty in an MSIP is represented by a finite scenario tree $\mathcal{T}$ over $T$ stages. Each node $n \in \mathcal{T}$ corresponds to a particular realization of uncertainty at a given stage. The children $C(n)$ of the node represent the possible outcomes in the next stage, and the probability of transitioning from parent node $n$ to child node $m$ is denoted by $q_{nm}$. At every node $n$, a $d_n$-dimensional state variable vector $x_n \in X_n$ links the decisions at that node to decisions at its descendants, while a local variable vector $y_n \in Y_n$ captures node-specific decisions. The sets $X_n$ and $Y_n$ may include integrality restrictions. The feasible decisions at node $n$, given the state $x_{a(n)}$ from the ancestor node $a(n)$, are described by a polyhedron $H_n(x_{a(n)})$, which depends on $x_{a(n)}$ through the linking constraints that connect decisions across stages; by contrast, the sets $X_n$ and $Y_n$ only depend on the local variable vector $y_n$. The objective is to minimize the total expected cost across all stages, yielding a recursive formulation in which each node's value function $Q_n(x_{a(n)})$ captures the minimum expected cost from that node onward:
\begin{subequations} \label{prob:sub}
\begin{align}
Q_n(x_{a(n)})
= \min_{x_n, y_n, z_n}\
& f_n(x_n, y_n)
+ \sum_{m \in C(n)} q_{nm}\, Q_m(x_n), \label{eq:obj} \\
& (x_n, y_n) \in H_n(z_n) \cap (X_n \times Y_n), \label{eq:feas} \\
& z_n = x_{a(n)}, \quad z_n \in Z_{a(n)} \supseteq X_{a(n)}. \label{eq:copy}
\end{align}
\end{subequations}
The objective function $f_n(\cdot)$ is assumed to be linear, and the initial state $x_0$ is given a priori. For nodes at the final stage, the set of children is empty, and the expected future cost is $0$. The function
$Q_n : \mathbb{R}^{d_{a(n)}} \to \mathbb{R} \cup \{+\infty\}$ takes value $\infty$
whenever the feasible set $H_n(x_{a(n)})$ is empty, and we denote by
$\operatorname{dom}(Q_n)$ the set of states $x_{a(n)}$ for which the subproblem $Q_n$ is feasible. For a set $S \subseteq \mathbb{R}^{d_{a(n)}}$, we write $\epi_S(Q_n) := \{(x_{a(n)}, \theta_n) \in S \times \mathbb{R} : \theta_n \geq Q_n(x_{a(n)})\}$ for the epigraph of $Q_n$ restricted to $S$. The copy constraint \eqref{eq:copy} with relaxed domain $Z_{a(n)} \supseteq X_{a(n)}$ plays an important role in cut generation. By dualizing the constraint, decomposition methods derive cuts of the form $\theta_n \geq g_n(x_{a(n)})$ that are valid for $\epi_{Z_{a(n)}}(Q_n)$ and are iteratively added to tighten the approximation of the value function. We denote by $\underline{Q}_n(\cdot)$ the current lower-bounding approximation of the value function $Q_n$, obtained by replacing each child value function $Q_m(x_n)$ in \eqref{prob:sub} with a cutting-plane approximation.

Although MSIPs admit a deterministic equivalent mixed-integer programming (MIP) formulation that bundles all stages and scenarios together, the number of variables and constraints grows exponentially with the scenario tree size, quickly overwhelming general-purpose MIP solvers. This motivates the use of Benders-type decomposition methods (\citealt{van1969shaped, birge1985decomposition, pereira1991multi}), which iteratively approximate the value functions $Q_n(\cdot)$ using cutting planes obtained from smaller scenario-subproblems. For problems with purely continuous recourse, classical Benders cuts (\citealt{benders1962partitioning}) derived from the linear programming (LP) dual of subproblems are sufficient for exact representation of the value function. However, when decision variables are integer, the value function is typically nonlinear and nonconvex (\citealt{blair1982value}), and LP-based cuts provide only a weak approximation. Specialized cuts have been proposed to address this, but each has limitations. Methods such as integer L-shaped cuts (\citealt{laporte1993integer}), disjunctive programming cuts (\citealt{sen2005c, sen2006decomposition}), parametric Gomory cuts (\citealt{gade2014decomposition, zhang2014finitely}), and combinatorial Benders cuts (\citealt{codato2004combinatorial}) assume purely binary or integer state variables and do not extend to general mixed-integer settings. Others rely on special problem structure absent in general MSIPs: logic-based Benders cuts (\citealt{hooker2023logic}) require problem-specific logical structure, MIDAS (\citealt{philpott2020midas}) requires monotone cost-to-go functions, and \citealt{ahmed2022stochastic} and \citealt{fullner2024lipschitz} assume Lipschitz continuity of the value function. Still others face scalability limitations. For example, scaled cuts (\citealt{van2024converging}) converge to the convex envelope of the expected recourse function in theory, but cannot be computed via scenario decomposition, limiting their scalability in the multistage setting (\citealt{romeijnders2024benders}). To the best of our knowledge, no available implementation of these methods has been shown to scale to general MSIPs with mixed-integer state variables.

Ensuring convergence for general MSIPs requires cuts that are \emph{tight} even when state variables are mixed-integer. A cut $\theta_n \geq g_n(x_{a(n)})$ is \emph{tight} at an incumbent $\hat{x}_{a(n)}$ if it holds at equality at that point, i.e., $g_n(\hat{x}_{a(n)}) = \un{Q}_n(\hat{x}_{a(n)})$. \citealt{zou2019stochastic} propose Lagrangian cuts derived by relaxing the copy constraints \eqref{eq:copy}. When state variables are binary, Lagrangian cuts obtained from optimal solutions of the Lagrangian dual are guaranteed to be tight, yielding finite convergence of the decomposition method. For non-binary state variables, however, \citealt{zou2019stochastic} resort to encoding the continuous and general integer variables with auxiliary binary variables, which can dramatically expand the state space and make the Lagrangian dual computationally prohibitive. \citealt{deng2024relu} address this limitation by reformulating the copy constraints using the nonlinear rectified linear unit (ReLU) functions, which partition the state space around the incumbent and keep the number of dual variables bounded by twice the state dimension. The resulting \emph{ReLU dual} admits strong duality under standard assumptions, so cuts generated from its optimal solutions are guaranteed to be tight at the incumbent even with mixed-integer state variables. This property ensures asymptotic convergence for general MSIPs, making ReLU cuts an important tool for multistage problems.

Despite this convergence guarantee, the ReLU dual typically admits multiple optimal solutions, and more generally, dual degeneracy can arise in other dual formulations, including Lagrangian and LP duals. Each optimal dual solution produces a cut that is tight at the current incumbent but may provide a poor approximation of the value function at other points in the state space. For example, \citealt{bansal2026integer} show that the coefficients of integer L-shaped cuts, known to give weak global approximations, are optimal solutions to the Lagrangian dual and since the ReLU dual coincides with the Lagrangian dual for binary state variables, the same holds for the ReLU dual.

Two main strategies mitigate this dual degeneracy problem across LP, Lagrangian, and ReLU duals. The first is \emph{normalization}: additional constraints on the dual variables select solutions that yield stronger cuts (\citealt{fischetti2010note, chen2022generating, fullnernew}). \citealt{bansal2026normalization} extend this framework to the ReLU dual, proving that normalized ReLU cuts are tight and Pareto-optimal in the original state space. The second is \emph{regularization}, where an augmented objective based on a relative interior point of the current relaxation, referred to as a \emph{core point}, selects an optimal dual solution that yields a Pareto-optimal cut (\citealt{magnanti1981accelerating, yangyang2025}). However, \citealt{bansal2026normalization} show that normalization is more flexible, as any cut obtained via regularization can also be obtained by normalization, but the converse does not hold. A complementary acceleration strategy is the \emph{alternating cut criterion} of \citealt{angulo2016improving} where cheaper Benders cuts from the LP relaxation are added whenever they cut off the incumbent, and the expensive ReLU or Lagrangian dual is solved only when they fail. This lowers the per-iteration cost, and the computational study in \citealt{bansal2026integer} shows that Benders cuts can also help mitigate dual degeneracy.

Several software packages exist for solving multistage stochastic programs. SDDP.jl (\citealt{dowson2021sddp}), the most widely used open-source library, implements both SDDP and the SDDiP algorithm of \citealt{zou2019stochastic}, but its integer support is limited to binary state variables via binary expansion. Among the alternatives are msppy (\citealt{ding2019python}), StOpt (\citealt{gevret2018stochastic}), StochDynamicProgramming.jl (\citealt{StochDynamicProgramming.jl}), and DynamicSDDP.jl (\citealt{fullnernew}); none of them support cutting planes that ensure convergence with general mixed-integer state variables. Cut selection strategies designed to improve approximation quality, such as Pareto-optimal cuts (\citealt{magnanti1981accelerating}), the MIP-based unified framework (\citealt{fischetti2010note}), and deepest cuts (\citealt{hosseini2025deepest}), have been developed and implemented for LP-based Benders subproblems or for specific problem classes such as facility location (\citealt{cordeau2019benders}) or energy planning (\citealt{psr-sddp}), but not for more general families of duals for MIP subproblems, such as ReLU duals, in general multistage stochastic programs.

In this paper, we present an open-source software package that implements a Benders-type decomposition framework for MSIPs, incorporating a recently proposed family of cuts that ensure convergence with general mixed-integer state variables, along with other cut-strengthening strategies. Specifically, the package builds on top of the SDDP.jl package (\citealt{dowson2021sddp}) and incorporates the following components:
\begin{enumerate}
    \item \textbf{ReLU cuts}: We implement the ReLU cuts (\citealt{deng2024relu}), which provide convergence guarantees for MSIPs with general mixed-integer state variables.
    \item \textbf{Cut quality and efficiency}: We implement three strategies to improve cut strength and reduce per-iteration cost: (i) the normalization framework (\citealt{bansal2026normalization}, \citealt{fullnernew}) for selecting dual solutions that lead to stronger cuts, (ii) the regularization approach of \citealt{yangyang2025} as an alternative for resolving dual degeneracy, and (iii) the alternating cut criterion of \citealt{angulo2016improving}, which alternates between tight and computationally efficient cuts to accelerate convergence.
    \item \textbf{Extensive experiments}: We benchmark all implemented methods on four classes of real-world MSIPs using publicly available instances---capacitated lot-sizing, generation expansion planning, airline revenue management, and portfolio optimization---providing insights into how problem structure affects the practical performance of each strategy.
\end{enumerate}

Together, these components make the package the first open-source implementation with (asymptotic) convergence guarantees for general MSIPs with mixed-integer state variables. It provides both the theoretical tools needed for convergence and the acceleration strategies needed for competitive computational performance. The package is designed to be versatile, allowing researchers and practitioners to configure key algorithmic parameters, such as the cut family (Benders, strengthened Benders, Lagrangian, or ReLU), the order in which cuts are added (alternating cut criterion), and the cut-strengthening approach (normalization or regularization).

The remainder of this paper is organized as follows. Section~\ref{sec:implementation} provides background on the methods implemented, covering the Lagrangian dual, the ReLU dual, normalization, regularization, and the alternating cut criterion. Section~\ref{sec:implementation_details} describes the implementation details, including the assumptions required by the algorithm, the construction of core points, and how to pass input in the software package. Section~\ref{sec:comp} presents the computational study. Section~\ref{sec:conclusion} concludes with a summary of findings and takeaways from our implementation and experiments. Appendix~\ref{sec:appendix} provides the formulations for the four problem classes considered in the experiments.

\section{Dual Formulations and Cut Strengthening} \label{sec:implementation}

\subsection{Lagrangian Dual and Lagrangian Cuts} \label{sec:lagrangian_dual}

Consider the subproblem \eqref{prob:sub} and the copy constraint $z_n = x_{a(n)}$ in \eqref{eq:copy}. The Lagrangian dual is obtained by relaxing this constraint with a multiplier $\pi_n \in \mathbb{R}^{d_{a(n)}}$, yielding the Lagrangian relaxation:
\begin{align} \label{eq:lagrn_relax}
    \mathcal{L}^O_n(\pi_n;\, \hat{x}_{a(n)}) := \min_{z_n \in Z_{a(n)}} \left[ \underline{Q}_n(z_n) + \pi_n^\top z_n - \pi_n^\top \hat{x}_{a(n)} \right],
\end{align}
where the superscript $O$ denotes the original Lagrangian based on the standard copy constraint \eqref{eq:copy}. The corresponding Lagrangian dual problem is:
\begin{align} \label{eq:lagrn_dual}
    \max_{\pi_n \in \mathbb{R}^{d_{a(n)}}}  \mathcal{L}^O_n(\pi_n;\, \hat{x}_{a(n)}).
\end{align}
For any dual multiplier $\pi_n$, the following Lagrangian cut is valid for $\epi_{Z_{a(n)}}(Q_n)$:
\begin{align*}
    \theta_n \geq \mathcal{L}^O_n(\pi_n;\, \hat{x}_{a(n)}) + \pi_n^\top (\hat{x}_{a(n)} - x_{a(n)}).
\end{align*}
\citealt{zou2019stochastic} prove that when the state variables are binary, strong duality holds for the Lagrangian dual \eqref{eq:lagrn_dual}, i.e., the optimal value of \eqref{eq:lagrn_dual} equals $\underline{Q}_n(\hat{x}_{a(n)})$. As a result, cuts generated from optimal solutions of \eqref{eq:lagrn_dual} are tight at the incumbent, guaranteeing finite convergence of the decomposition method. However, when the state variables are general integer or mixed-integer, strong duality need not hold for the Lagrangian dual based on \eqref{eq:copy}, and the resulting cuts are not guaranteed to be tight. This limitation motivates the ReLU-based reformulation described next.

\subsection{ReLU Dual and ReLU Cuts} \label{sec:relu_dual}

To address the challenge of the Lagrangian dual with general mixed-integer state variables, \citealt{deng2024relu} replace the copy constraint \eqref{eq:copy} with an equivalent ReLU-based reformulation:
\begin{align} \label{eq:relu_copy}
    (z_{nk} - x_{a(n)k})^+ = 0, \quad (z_{nk} - x_{a(n)k})^- = 0, \quad \forall k \in [d_{a(n)}],
\end{align}
where $(x)^+ := \max\{x, 0\}$ and $(x)^- := \max\{-x, 0\}$. These ReLU functions partition the state space around the incumbent, with each partition reflecting the space to either left or right of the incumbent along some dimension $k \in [d_{a(n)}]$. Dualizing these ReLU-based copy constraints, rather than the equivalent linear constraints \eqref{eq:copy}, allows strong duality to hold for the resulting dual problem.

Dualizing the constraints \eqref{eq:relu_copy} yields the following Lagrangian relaxation, parameterized by dual multipliers $\pi_n^+, \pi_n^- \in \mathbb{R}^{d_{a(n)}}$:
\begin{align*}
    \mathcal{L}^R_n(\pi_n^+, \pi_n^-; \hat{x}_{a(n)}) := \min_{z_n \in Z_{a(n)}} \left[ \underline{Q}_n(z_n) + \sum_{k=1}^{d_{a(n)}} \pi_{nk}^+ (z_{nk} - \hat{x}_{a(n),k})^+ + \sum_{k=1}^{d_{a(n)}} \pi_{nk}^- (z_{nk} - \hat{x}_{a(n),k})^- \right],
\end{align*}
where $g_n(x_{a(n)})$ in the cut expression is obtained by solving the corresponding ReLU dual problem:
\begin{align} \label{eq:relu_dual}
    \max_{\pi_n^+, \pi_n^- \in \mathbb{R}^{d_{a(n)}}} \mathcal{L}^R_n(\pi_n^+, \pi_n^-; \hat{x}_{a(n)}).
\end{align}
For any dual multipliers $\pi_n^+, \pi_n^-$, the following ReLU cut is valid for $\epi_{Z_{a(n)}}(Q_n)$:
\begin{align} \label{eq:relu_cut}
    \theta_n \geq \mathcal{L}^R_n(\pi_n^+, \pi_n^-; \hat{x}_{a(n)}) - \sum_{k=1}^{d_{a(n)}} \pi_{nk}^+ (x_{a(n),k} - \hat{x}_{a(n),k})^+ - \sum_{k=1}^{d_{a(n)}} \pi_{nk}^- (x_{a(n),k} - \hat{x}_{a(n),k})^-.
\end{align}
\citealt{deng2024relu} prove that strong duality holds for the ReLU dual \eqref{eq:relu_dual}, i.e., the optimal value of \eqref{eq:relu_dual} equals $\underline{Q}_n(\hat{x}_{a(n)})$, so ReLU cuts generated from optimal dual solutions are tight at the incumbent even when state variables are mixed-integer, a property that ensures asymptotic convergence for general MSIPs.

Strong duality alone does not guarantee strong cuts, as dual degeneracy can yield different optimal solutions to \eqref{eq:relu_dual} that yield cuts with widely varying global approximation quality. This motivates the cut-strengthening strategies described below.

\subsection{Normalized Dual} \label{sec:normalization}

Normalization addresses dual degeneracy by imposing additional constraints on the dual variables to select solutions that yield stronger cuts. For LP duals, \citealt{fischetti2010note} discuss a choice of normalization constraint that yields cuts corresponding to a minimal infeasible subsystem of the scenario subproblem. \citealt{chen2022generating} and \citealt{fullnernew} extend this idea to the Lagrangian dual (as in Section~\ref{sec:lagrangian_dual}) of MIP subproblems, showing that appropriately chosen normalization coefficients lead to cuts that are Pareto-optimal. We present normalization here in the context of the ReLU dual, but the concept applies analogously to other dual formulations.

\citealt{bansal2026normalization} develop a normalization framework for the ReLU dual. A key step in developing this framework is to reinterpret the ReLU cut \eqref{eq:relu_cut} as a standard Lagrangian cut in a lifted space. Specifically, via the change of variables $w_{nk}^+ = (z_{nk} - \hat{x}_{a(n),k})^+$ and $w_{nk}^- = (z_{nk} - \hat{x}_{a(n),k})^-$, the ReLU cut at incumbent $\hat{x}_{a(n)}$ in the original space corresponds to a Lagrangian cut at the origin $(\mathbf{0}, \mathbf{0})$ in the lifted space of $(w_n^+, w_n^-)$. The lifted domain is defined as the image of $Z_{a(n)}$ under this change of variables:
\begin{align*}
    Z^{\mathrm{lift}}_{\hat{x}_{a(n)}} = \left\{ (w_n^+, w_n^-) : w_{nk}^+ = (z_{nk} - \hat{x}_{a(n),k})^+,\; w_{nk}^- = (z_{nk} - \hat{x}_{a(n),k})^-\; \forall k \in [d_{a(n)}],\; z_n \in Z_{a(n)} \right\}.
\end{align*}
Defining the lifted recourse function $\underline{Q}_n'(w_n^+, w_n^-;\, \hat{x}_{a(n)}) := \underline{Q}_n(\hat{x}_{a(n)} + w_n^+ - w_n^-)$, this correspondence enables applying the normalization procedure of \citealt{fullnernew} to the following feasibility version of the lifted subproblem at incumbent $(\mathbf{0}, \mathbf{0})$ and epigraph variable value $\hat{\theta}_n$:
\begin{subequations} \label{eq:feas_sub}
\begin{align}
    \min \quad & 0 \label{feas_sub:obj} \\
    \text{s.t.} \quad & \underline{Q}_n'(w_n^+, w_n^-;\, \hat{x}_{a(n)}) \leq \hat{\theta}_n, \quad (w_n^+, w_n^-) = (\mathbf{0}, \mathbf{0}), \label{eq:feas_theta_copy} \\
    & (w_n^+, w_n^-) \in Z^{\mathrm{lift}}_{\hat{x}_{a(n)}}. \label{eq:feas_lifted}
\end{align}
\end{subequations}
The feasibility subproblem \eqref{eq:feas_sub} captures whether the incumbent $(\mathbf{0}, \mathbf{0})$ (corresponding to $\hat{x}_{a(n)}$ in the original space) is feasible in the epigraph of $\underline{Q}_n'$ at cost $\hat{\theta}_n$ and with respect to the constraint in \eqref{eq:feas_lifted}.

Introducing multipliers $\pi_{n0} \geq 0$ for the epigraph constraint and $(\pi_n^+, \pi_n^-)$ for the copy constraint in \eqref{eq:feas_theta_copy}, respectively, and normalization coefficients $(u_n^+, u_n^-, u_{n0})$, the normalized dual is:
\begin{subequations}\label{eq:norm_dual}
\begin{align}
    \max_{\pi_n^+, \pi_n^-, \pi_{n0}} \quad & \mathcal{L}_n(\pi_n^+, \pi_n^-, \pi_{n0}; \hat{x}_{a(n)}) - \pi_{n0} \hat{\theta}_n \label{eq:norm_dual_obj} \\
    \text{s.t.} \quad & (u_n^+)^\top \pi_n^+ + (u_n^-)^\top \pi_n^- + u_{n0} \pi_{n0} \leq 1, \quad \pi_{n0} \geq 0, \label{eq:norm_dual_constr}
\end{align}
\end{subequations}
where $\mathcal{L}_n(\pi_n^+, \pi_n^-, \pi_{n0}; \hat{x}_{a(n)})$ is the Lagrangian relaxation of the feasibility subproblem \eqref{eq:feas_sub} obtained by dualizing constraints \eqref{eq:feas_theta_copy}. Based on an optimal solution $(\hat{\pi}_n^+, \hat{\pi}_n^-, \hat{\pi}_{n0})$ of \eqref{eq:norm_dual}, the normalized ReLU cut is:
\begin{align*}
    \hat{\pi}_{n0}\, \theta_n \geq \mathcal{L}_n(\hat{\pi}_n^+, \hat{\pi}_n^-, \hat{\pi}_{n0}; \hat{x}_{a(n)}) - \sum_{k=1}^{d_{a(n)}} \hat{\pi}_{nk}^+ (x_{a(n),k} - \hat{x}_{a(n),k})^+ - \sum_{k=1}^{d_{a(n)}} \hat{\pi}_{nk}^- (x_{a(n),k} - \hat{x}_{a(n),k})^-.
\end{align*}
The strength of a normalized cut depends on the normalization coefficients $(u_n^+, u_n^-, u_{n0})$. \citealt{bansal2026normalization} define Pareto-optimality for nonlinear ReLU cuts over the reference set $\mathcal{H}_n$, the tightest epigraph obtained by intersecting all valid ReLU cuts at a given incumbent, and show that normalization using the right choice of coefficients $(u_n^+, u_n^-, u_{n0})$ yields a Pareto-optimal cut over this set. As discussed next, these coefficients are closely related to the regularized objective coefficients in the regularized dual, revealing an important connection between the normalized and regularized dual formulations.

\subsection{Regularized Dual} \label{sec:regularization}

\citealt{yangyang2025} propose an alternative strategy for selecting optimal solutions to the dual problem that lead to strong cuts. Although we present regularization here for the ReLU dual \eqref{eq:relu_dual}, the same idea applies to other dual formulations; for instance, replacing $\mathcal{L}^R_n$ with the original Lagrangian relaxation $\mathcal{L}^O_n$ from \eqref{eq:lagrn_relax} yields a regularized Lagrangian dual. The idea is to characterize the set of all optimal solutions $\Pi_n(\hat{x}_{a(n)})$ of the dual and optimize a regularized objective over this set. Given a core point $(\tilde{u}_n^+, \tilde{u}_n^-)$ in relative interior of $\conv(Z^{\mathrm{lift}}_{\hat{x}_{a(n)}})$, the regularized dual is given by:
\begin{align} \label{eq:reg_dual}
    \max_{\pi_n^+, \pi_n^- \in \Pi_n(\hat{x}_{a(n)})} \left[ \mathcal{L}^R_n(\pi_n^+, \pi_n^-; \hat{x}_{a(n)}) - (\pi_n^+)^\top \tilde{u}_n^+ - (\pi_n^-)^\top \tilde{u}_n^- \right].
\end{align}
The optimal set $\Pi_n(\hat{x}_{a(n)})$ is characterized implicitly by the constraint $\mathcal{L}^R_n(\pi_n^+, \pi_n^-; \hat{x}_{a(n)}) \geq \underline{Q}_n(\hat{x}_{a(n)}) - \epsilon$ for a small $\epsilon > 0$, and approximated iteratively via the level-bundle method since $\mathcal{L}^R_n$ is concave and piecewise linear.

\citealt{yangyang2025} prove that the cut resulting from \eqref{eq:reg_dual} is both tight at the incumbent and Pareto-optimal in the lifted space $(w_n^+, w_n^-)$, specifically for the convex envelope of the lifted recourse function $\underline{Q}'_n$ over $\conv(Z^{\mathrm{lift}}_{\hat{x}_{a(n)}})$. Because the dual variables are constrained to lie in $\Pi_n(\hat{x}_{a(n)})$, the regularized dual always produces cuts that are tight at the incumbent.

The objective in \eqref{eq:reg_dual} can also be interpreted as follows. While enforcing tightness at the current incumbent through
$\pi_n^+, \pi_n^- \in \Pi_n(\hat{x}_{a(n)})$, it also seeks a valid cut that is as tight as possible at another interior point of the domain, namely
$(\ti{u}_n^+, \ti{u}_n^-)$. \citealt{fullnernew} discuss that the normalization coefficients admit a similar interpretation. In particular, the corresponding core point is given by $((\mbf{0},\mbf{0}),\hat{\theta}_n) + \eta_n^* (u_n^+,u_n^-,u_{n0})$, where $\eta_n^*$ is an optimal solution of the normalized dual \eqref{eq:norm_dual}. Owing to the one-to-one relationship between the core points and the normalization coefficients, we refer to the normalization coefficients as core points throughout the rest of the paper.

\citealt{bansal2026normalization} establish another important relation between normalization and regularization: any cut obtained by the regularized dual \eqref{eq:reg_dual} can also be obtained by the normalized dual \eqref{eq:norm_dual} with an appropriate choice of normalization coefficients. The converse, however, does not hold. Regularization is restricted to optimal solutions of the ReLU dual, i.e., to tight cuts, whereas normalization can also generate Pareto-optimal cuts that are not tight at the incumbent. Such cuts may provide stronger global approximations of the value function by better covering the state space at points other than the incumbent, yet still cut off the incumbent. This flexibility gives normalization an advantage over regularization, which we depict in our computational experiments.

\subsection{Alternating Cut Criterion} \label{sec:alternating}

Solving the ReLU dual, whether in its standard, normalized, or regularized form, requires multiple evaluations of the mixed-integer subproblem and can be computationally expensive. The alternating cut criterion of \citealt{angulo2016improving} addresses this by deferring to a cheaper Benders cut derived from the LP relaxation whenever it can separate the incumbent, and using the expensive ReLU or Lagrangian cuts only when necessary. Since the LP relaxation is far less expensive to solve, this can reduce per-iteration cost.

The alternating criterion was originally proposed in the two-stage setting, using Benders cuts and integer L-shaped cuts, with the primary goal of obtaining computational speedups, but in our setting, we view the idea more broadly. One may alternate between any two families of cuts, provided that one family guarantees convergence and the other can be generated cheaply to cut off the incumbent or improve the relaxation.

This broader interpretation has an additional benefit beyond efficiency. The cheaper cuts need not be tight at the incumbent, but they are not necessarily dominated by tight cuts. Tight cuts accurately approximate the value function at the incumbent, yet they may provide weak lower bounds elsewhere in the state space. By contrast, LP-based Benders cuts, and potentially other inexpensive cut families, can yield stronger approximations at some points away from the incumbent. Thus, alternating cuts can exploit this complementarity. Tight cuts ensure convergence, while cheaper cuts can improve the global approximation of the value function at low computational cost.

\section{Implementation Details} \label{sec:implementation_details}

This section covers the implementation details of the package. We first state the assumptions required for the correctness and convergence of the algorithm. We then describe how the core points used by the normalized and regularized duals are computed in practice. Finally, we explain how to pass input to the algorithm to configure the hyperparameters.

\subsection{Assumptions} \label{sec:assumptions}

We make two standard assumptions: first, the set $X_n$ is assumed to be compact. Under this assumption, the set $X_n$ can be shifted so that all state variables satisfy $x_{nk} \in [0, B_k], k \in [d_{a(n)}]$ for finite bounds $B_k$. Second, we assume $\dom(Q_n) = Z_{a(n)}$ is of the form $\prod_{k} [0, B_k]$. Together with the domain relaxation in \eqref{eq:copy}, this guarantees relatively complete recourse, that is, for every feasible first-stage decision, there exists at least one feasible recourse action in subsequent stages regardless of the realized scenario, a standard assumption in stochastic programming. This condition can always be enforced by adding continuous variables to $H_n$ and penalizing them in the objective. The set $Z_{a(n)}$ is typically taken as a superset of $X_{a(n)}$ to reduce computational effort in solving the dual.

\subsection{Core Point Computation} \label{sec:core_point}

As discussed in Section \ref{sec:regularization}, since there is a one-to-one correspondence between normalization coefficients $(u_n^+, u_n^-, u_{n0})$ and the regularization core point $(\ti{u}_n^+, \ti{u}_n^-)$, we use the term core point for both, and we now describe how the core points are constructed in practice. Under the assumptions discussed in previous subsection, $Z_{a(n)}$ is of the form $\prod_k [0, B_k]$, so the lifted domain $Z^{\mathrm{lift}}_{\hat{x}_{a(n)}}$ decomposes dimension-wise. Specifically, Proposition~3 in \citealt{bansal2026normalization} establishes that $\conv(Z^{\mathrm{lift}}_{\hat{x}_{a(n)}}) = \prod_k \conv(Z^{\mathrm{lift}}_{\hat{x}_{a(n)},k})$, where each factor is the triangle $T_k := \conv\!\bigl\{(0,0),\;(B_k - \hat{x}_{a(n),k},\, 0),\;(0,\, \hat{x}_{a(n),k})\bigr\}$.
The origin $(\mbf{0}, \mbf{0})$ in the lifted space corresponds to the incumbent $z_{nk} = \hat{x}_{a(n),k}$ in the original space, the vertices $(B_k - \hat{x}_{a(n),k}, 0)$ and $(0, \hat{x}_{a(n),k})$ correspond to $z_{nk} = B_k$ and $z_{nk} = 0$, respectively.

Proposition~3 of \citealt{bansal2026normalization} requires $(u_{nk}^+, u_{nk}^-) \in \relint(T_k)$ for each $k \in [d_{a(n)}]$, together with $u_{n0} > \underline{Q}'_n(u_{n}^+, u_{n}^-; \hat{x}_{a(n)}) - \hat{\theta}_n$, to obtain a Pareto-optimal cut; Proposition~4 further shows that, for an interior incumbent coordinate $0 < \hat{x}_{a(n),k} < B_k$, we may impose $u_{nk}^+ = u_{nk}^-$. This symmetry is useful because it keeps $\hat{x}_{a(n)} + u_n^+ - u_n^-$ equal to the incumbent, so $u_{n0}$ can be approximated from the already available value $\underline{Q}_n(\hat{x}_{a(n)})$ without an additional subproblem evaluation. Thus, writing $u_{nk}^+ = u_{nk}^- = c$, the diagonal point $(c,c)$ lies in $\relint(T_k)$ if and only if
\begin{align*}
    0 < c < e_k, \qquad \text{where} \quad e_k := \frac{\hat{x}_{a(n),k}(B_k - \hat{x}_{a(n),k})}{B_k}.
\end{align*}
The scalar $e_k$ is the largest admissible symmetric value; it is maximized at the midpoint $\hat{x}_{a(n),k} = B_k/2$ and vanishes when $\hat{x}_{a(n),k}$ is at either bound.

We implement two core-point families in our package. They handle three cases depending on whether the state is at the upper bound ($\hat{x}_{a(n),k} = B_k$), at the lower bound ($\hat{x}_{a(n),k} = 0$), or in the interior ($0 < \hat{x}_{a(n),k} < B_k$). At the boundary, $e_k = 0$, so the non-zero component is instead set proportional to the distance to the opposite bound. \emph{CoreEps}, parameterized by $(\varepsilon_{\mathrm{bnd}}, \varepsilon_{\mathrm{int}})$ with $0 < \varepsilon_{\mathrm{bnd}}, \varepsilon_{\mathrm{int}} \ll 1$, places the core point near the origin of $T_k$:
\begin{align*}
    u_{nk}^+ &= \begin{cases}
        0 & \hat{x}_{a(n),k} = B_k, \\
        \varepsilon_{\mathrm{bnd}}\,(B_k - \hat{x}_{a(n),k}) & \hat{x}_{a(n),k} = 0, \\
        \varepsilon_{\mathrm{int}}\, e_k & \text{otherwise},
    \end{cases} \qquad &
    u_{nk}^- &= \begin{cases}
        \varepsilon_{\mathrm{bnd}}\,\hat{x}_{a(n),k} & \hat{x}_{a(n),k} = B_k, \\
        0 & \hat{x}_{a(n),k} = 0, \\
        \varepsilon_{\mathrm{int}}\, e_k & \text{otherwise}.
    \end{cases}
\end{align*}
\emph{CoreMid}, parameterized by $(\sigma_{\mathrm{bnd}}, \sigma_{\mathrm{int}})$ with $0 < \sigma_{\mathrm{bnd}}, \sigma_{\mathrm{int}} < 1$, uses the same structure ($\epsilon$ replaced with $\sigma$) but with larger scaling factors, placing the core point deeper inside $T_k$. At $\sigma_{\mathrm{int}} = 0.5$, the interior component equals $e_k/2$, placing the core point at the midpoint of the admissible range along the diagonal.

For normalization, the pair $(u_{nk}^+, u_{nk}^-)$ computed by either family is used in the normalization constraint \eqref{eq:norm_dual_constr}. The normalization constraint also includes an objective component $u_{n0} = \underline{Q}_n(\hat{x}_{a(n)}) - \hat{\theta}_n + \delta$, with $\delta > 0$, where $\underline{Q}_n(\hat{x}_{a(n)})$ is the value returned by the current subproblem solve, $\hat{\theta}_n$ is the current epigraph variable value at the parent node, and $\delta$ is an additional hyperparameter (default $\delta = 1$). This satisfies the condition $u_{n0} > \underline{Q}_n(\hat{x}_{a(n)}) - \hat{\theta}_n$ from Proposition 3 of \citealt{bansal2026normalization}. For regularization, the same $(u_{nk}^+, u_{nk}^-)$ serve as $(\ti{u}_{nk}^+, \ti{u}_{nk}^-)$ in \eqref{eq:reg_dual}; no objective component is needed.

The parameter $\varepsilon_{\mathrm{int}}$ (resp.\ $\sigma_{\mathrm{int}}$) controls the fraction of $e_k$ used as the core-point coordinate. A value close to $0$ places the core point near the origin of $Z^{\mathrm{lift}}_{\hat{x}_{a(n)}}$ and tends to produce cuts that are tighter at the current incumbent, while a value close to $1$ places it deeper in the interior and tends to yield stronger global approximations of the value function. Because neither family dominates the other a priori, the computational study evaluates both and selects the best hyperparameter setting on a held-out tune instance; see Section~\ref{sec:comp}.

\subsection{Configuring the Algorithm} \label{sec:configuration}

The methods described in this paper are implemented as a Julia package, \texttt{SDDmiP.jl}, available at \url{https://github.com/akulbansal5/SDDmiP}, built on top of the SDDP.jl package (\citealt{dowson2021sddp}); the interface follows the same variable names as the original code. The cut family used in the method is selected by passing a \texttt{duality\_handler} keyword argument to: \linebreak
\texttt{SDDmiP.train(model; duality\_handler = handler)}.

The original package already provides support for Benders, strengthened Benders, and Lagrangian cuts using duality handlers \texttt{ContinuousConicDuality()}, \texttt{StrengthenedConicDuality()}, and \linebreak \texttt{LagrangianDuality()}. We add support for ReLU cuts and cut-strengthening methods using \linebreak \texttt{LagrnNormDuality()}, \texttt{LagrnRegDuality()}, \texttt{ReLUNormDuality()} and \texttt{ReLURegDuality()}.

The normalized and regularized handlers accept an optional \texttt{core\_type} keyword argument that specifies the core point used in the dual, as described in Section~\ref{sec:core_point}. Passing \texttt{CoreEps(tol, $\varepsilon_{\mathrm{bnd}}$, $\varepsilon_{\mathrm{int}}$, $\delta$)} places the core point near the origin of each triangle $T_k$, while \texttt{CoreMid(tol, $\sigma_{\mathrm{bnd}}$, $\sigma_{\mathrm{int}}$, $\delta$)} places it deeper in the interior. The first argument \texttt{tol} is a numerical tolerance that classifies each incumbent coordinate as at the upper bound ($\hat{x}_{a(n),k} \ge B_k - \texttt{tol}$), at the lower bound ($\hat{x}_{a(n),k} \le \texttt{tol}$), or in the interior; this avoids treating floating-point incumbents near a bound as interior points, where the symmetric core-point formulas of Section~\ref{sec:core_point} would not apply. The remaining arguments set the boundary and interior scaling factors ($\varepsilon_{\mathrm{bnd}}$, $\varepsilon_{\mathrm{int}}$ or $\sigma_{\mathrm{bnd}}$, $\sigma_{\mathrm{int}}$) and the objective offset $\delta$ used in the definition of $u_{n0}$.

\texttt{AlternativeBiDuality} implements the alternating cut criterion described in Section~\ref{sec:alternating}. It accepts a vector of exactly two handlers and alternates between the corresponding cut families at each backward-pass node. The first handler is attempted first; the second is invoked only if the first handler's cut fails to separate the incumbent. Typically, the first handler is \texttt{ContinuousConicDuality()}, and the second is either the Lagrangian or ReLU handler, so that cheaper Benders cuts are used whenever they suffice, and more expensive tight cuts are reserved for when needed.

\section{Computational Study}\label{sec:comp}

The computational study is designed to answer three questions: (i) how do normalization and regularization compare in terms of cut strength, solution time, overall convergence, and the computational effort required to solve the dual problem; (ii) when does the alternating criterion, which alternates with Benders cuts, improve performance and (iii) how does the structure of the underlying problem affect these comparisons.

We implement all methods within the stochastic dual dynamic programming (SDDP) framework of \citealt{pereira1991multi}, built on top of the Julia package SDDP.jl (\citealt{dowson2021sddp}). In multistage problems, the scenario tree grows exponentially with the number of stages, making it impractical to enumerate all paths in each iteration as nested Benders decomposition does (\citealt{birge1985decomposition}). SDDP addresses this by sampling a small subset of scenario paths in the forward pass and computing cuts only along those sampled paths in the backward pass (\citealt{pereira1991multi}). This sampling makes each iteration computationally tractable but introduces randomness. The lower bound from the master problem is deterministic, but the upper bound must be estimated from a separate out-of-sample evaluation and is therefore only a statistical estimate of the true upper bound. We terminate each run using a bound-stalling stopping rule: if the relative improvement in the lower bound over the preceding 10 iterations falls below $1\%$, the algorithm stops, subject to a wall-clock time limit of $3600$ seconds. The normalized and regularized ReLU duals are solved using the level bundle method of \citealt{lemarechal1995new}. All experiments are conducted on an Intel Xeon Gold 5218R CPU with 8 cores and 64 GB of RAM, and the LP and MIP subproblems are solved using Gurobi 12.0.

The experiments cover four classes of multistage stochastic integer programs: (i) generation expansion planning (GEP), in which the decision maker selects which power generators to build at each stage and how to dispatch them to meet uncertain demand while minimizing construction, generation, and lost-demand costs (\citealt{jin2011modeling}); (ii) airline revenue management (ARM), a network revenue management problem in which the decision maker determines how to allocate seats across fare classes and routes at each stage to maximize expected revenue from stochastic demand (\citealt{moller2008airline}); (iii) capacitated lot-sizing (CLSP), a production planning problem in which setup and production quantities are determined at each stage to meet stochastic demand at minimum cost (\citealt{trigeiro1989capacitated}); and (iv) portfolio optimization, a financial planning problem in which binary investment decisions are made at each stage to maximize expected returns (\citealt{zou2019stochastic}). Detailed formulations for each problem class are provided in Appendix~\ref{sec:appendix}.

For each problem class, instances are defined by a tuple of structural parameters (e.g., the number of stages, scenarios per stage, and the number of products or assets). Within each parameter configuration, we generate one \emph{tune} instance and three \emph{test} instances using the data generation procedure as in the associated reference, each with a different random seed for both the problem data and the scenario paths sampled during the SDDP algorithm. We report average results across the three test instances, which guards against conclusions specific to a single data realization or a single realized sequence of sampled paths.

Both normalization and regularization require a core point that lies in the relative interior of the domain. We consider two families of core points: \emph{CoreEps} and \emph{CoreMid}. Each family has hyperparameters that control exactly how close to the boundary, or how far into the middle, the core point is placed. The construction of the core point is discussed in detail in Section \ref{sec:core_point}.

Since both methods are sensitive to the choice of the core point, and each may perform best with a different family and parameter setting, we tune these hyperparameters on the tune instance before evaluating on held-out test instances. For each method and each configuration of core point family and parameters, we run the full algorithm on the tune instance and select the configuration that achieves the highest lower bound within a $2\%$ relative tolerance, breaking ties by shortest solve time. This two-phase tune-then-test protocol ensures that reported test metrics reflect hyperparameters chosen without access to the test data. With a larger computational budget, a finer grid search could yield further improvements.

Each test table reported below has the following columns: $T$ and $R$ denote the number of stages and scenarios per stage; \emph{Type} indicates whether the alternating criterion is disabled (D) or enabled (A); \emph{Method} is normalization (Norm) or regularization (Reg); \emph{LB} or \emph{UB} denotes the deterministic lower or upper bound at termination, depending on whether the problem is a minimization or maximization problem; \emph{Time} is total wall-clock time in seconds and \emph{TL} in the Time column indicates that the run hit the 3600-second time limit; \emph{Iters} is the number of SDDP iterations; \emph{Dual-Iters} and \emph{Dual t} is the average number of level-bundle iterations and average time in seconds required each time the dual problem is solved; \emph{Prop} denotes the fraction of added cuts that are ReLU cuts (relevant only under the alternating criterion); \emph{Alt Gap} is the average LP--MIP subproblem gap $|\text{MIP} - \text{LP}|/|\text{MIP}|$ observed when the alternating criterion is used, where MIP and LP are the optimal values of the MIP and LP relaxation of the subproblem, respectively; and \emph{Gap} denotes the relative optimality gap $|\text{UB} - \text{LB}|/|\text{UB}|$, where UB is the upper bound and LB is the best lower bound. All metrics reported in the test tables are averaged over the three test instances. The following discussion explains how this information helps answer the three questions posed at the beginning of this section.

\subsection{Generation Expansion Planning}

\input{tables/paper/gep_test.tex}

Table~\ref{tab:gep_test} reports results for the GEP instances. Under the default approach (Type = D), normalization (Norm) reduces solve time (Time (s)) by 34\% on average relative to regularization (Reg). Because our stopping criterion may also terminate when the bounds stall, lower solve times alone do not necessarily indicate faster convergence. We therefore also examine the final relative gap (Gap\%). In most cases, the gap under normalization is at most 0.40 percentage points larger than that of regularization. Normalization also requires at most as many iterations (Iters) in five of the eight cases. Taken together, fewer iterations and comparable gaps suggest that the cuts generated by normalization are as strong as those generated by regularization under the default approach. The lower solve times further indicate that the normalized dual is computationally easier to solve, making normalized ReLU cuts cheaper to obtain. The speed advantage of normalization under the default criterion stems entirely from the computational effort required to solve the dual problem. Average iterations to solve the dual (Dual-Iters) for normalization range from $25$--$30$ across all instance configurations, versus $41$--$45$ for regularization. Average dual time (Dual t(s)) is uniformly lower for Norm ($0.46$--$0.76$s versus $0.50$--$1.26$s for Reg).

Under the alternating criterion (A), the pattern of solution times reverses. Regularization method stalls faster than normalization in every row: for example, $20$s versus $26$s at $(10,3)$, $20$s versus $38$s at $(10,6)$, $28$s versus $55$s at $(10,12)$, and $30$s versus $144$s at $(10,20)$. However, these faster solve times often come with larger final gaps when regularization is applied. For instance, the final gap is $4.50\%$ for regularization versus $2.57\%$ for normalization at $(10,3)$, $3.19\%$ versus $1.88\%$ at $(10,10)$, and $5.68\%$ versus $2.44\%$ at $(10,12)$. The number of iterations and the fraction of ReLU cuts used (Prop.) are nearly identical across the two methods. Average iterations to solve the dual and average dual solve time is also less for regularization than normalization. This suggests that, under the alternating criterion, regularized ReLU cuts are cheaper to obtain but weaker than normalized ReLU cuts.

More importantly, the alternating criterion yields roughly an order-of-magnitude improvement in solve time, while increasing the final gaps by at most about one percentage point relative to the default criterion. The fraction of ReLU cuts (Prop.) among all added cuts is below $4\%$, indicating that cheaper Benders cuts close most of the gap and substantially reduce computational time. This is consistent with the small LP--MIP subproblem gaps reported in the Alt  Gap column, which are below $3\%$ in most cases. For normalization, the speedups under the alternating criterion stem primarily from solving the ReLU dual much less often, as reflected in the Prop. column. For regularization, the speedups come from both solving the ReLU dual less often and from reductions in the number of level-bundle iterations and dual-solve times relative to the default criterion.

\subsection{Airline Revenue Management}

\input{tables/paper/arm_test.tex}

Table~\ref{tab:arm_test} reports results for the ARM instances. Since ARM is a maximization problem, the reported bound is an upper bound (UB) on the optimal revenue. Across all $(T,R)$ configurations, we observe that regularization under default approach (Reg-D) fails to improve on the trivial upper bound, producing final gaps above $82\%$ in every case. Moreover, each Reg-D run terminates in exactly $11$ iterations, indicating that the upper bound stalls in first few iterations under the $1\%$ improvement criterion. This behavior illustrates a classic weakness of regularization-based cuts. Regularization selects Pareto-optimal cuts, but it also limits the selection to cuts that are tight at the current incumbent. Such cuts can be weak elsewhere in the state space. A non-tight Pareto-optimal cut may provide a better global approximation even though it is not tight at the incumbent. In these instances, the large upper bounds and final gaps indicate that the regularized cuts are notably weak away from the incumbent, so that a trivial upper bound dominates the bound obtained from these cuts.

In contrast, normalized ReLU cuts (Norm) under the default criterion (D) substantially tighten the upper bound. The final gaps under Norm-D range from $14.16\%$ to $25.69\%$, compared with gaps above $82\%$ under Reg-D. Norm-D also requires at most as many iterations as Reg-D in most cases. Together, these results indicate that normalized cuts provide a much stronger approximation of the value function on these instances. However, this improvement comes at a high computational cost: every Norm-D run hits the time limit. The bottleneck is the normalized dual, whose average solve time ranges from $27$ to $64$ seconds on the ARM instances, compared with less than $1$ second on the GEP instances. Thus, while normalization avoids the stagnation observed under Reg-D, the ARM subproblems make the normalized dual substantially harder to solve.

Finally, the alternating criterion (A) substantially improves both solve times and final gaps. For normalization, these gains come from two sources. First, the ReLU dual is solved for only a small fraction of the added cuts: the Prop. values are low, ranging from $0.4\%$ to $6.5\%$. Second, when the normalized dual is solved, it is cheaper under the alternating criterion. Both the average number of level-bundle iterations (Dual-Iters) and the average dual solve time (Dual t) decrease relative to Norm-D. Under the alternating approach, Norm-A and Reg-A achieve nearly identical final gaps and iteration counts across all configurations, indicating that the two methods generate cuts of similar strength in this setting. However, Norm-A is substantially faster than Reg-A, showing that normalized ReLU cuts are much cheaper to obtain on the ARM instances. The remaining gaps of $3\%$--$6\%$ under Norm-A reflect both bound stalling and sampling noise in the upper bound estimate.

\subsection{Capacitated Lot-Sizing}

\input{tables/paper/clsp_test.tex}

Table~\ref{tab:clsp_test} reports results for the CLSP instances. Under the default criterion (D), normalization (Norm) is consistently faster than regularization (Reg): for example, $134$s versus $173$s at $(P,T,R) = (3,5,2)$, $725$s versus $1330$s at $(6,5,2)$, $344$s versus $713$s at $(3,5,4)$, and $2012$s versus $2828$s at $(6,4,4)$. This speed advantage is driven by the normalized dual being cheaper to solve: Dual-Iters range from $21$--$47$ for Norm-D versus $24$--$60$ for Reg-D, and average Dual t is $0.54$--$5.27$s for Norm versus $0.65$--$6.81$s for Reg. The final gaps and iteration counts are comparable in most configurations, suggesting cuts of similar strength. In a few cases, Reg-D achieves a smaller gap: at $(3,5,4)$, Reg-D reaches $3.95\%$ versus Norm-D's $4.70\%$, and at $(6,4,8)$, Reg-D reaches $7.81\%$ versus Norm-D's $8.60\%$. However, these improvements come at substantially higher solution time (roughly $2\times$ at $(3,5,4)$, and Reg-D hits the time limit at $(6,4,8)$). Overall, Norm-D is more computationally efficient, while Reg-D occasionally achieves a better gap but at higher cost.

Under the alternating criterion (A), the relative ordering between Norm and Reg is less consistent than under the default criterion. In most configurations, Norm-A and Reg-A produce comparable gaps and iteration counts, with Norm-A being somewhat faster (e.g.\ $1060$s versus $2113$s at $(6,5,2)$, and $3544$s versus TL at $(6,5,4)$). In a few configurations, however, Reg-A produces significantly larger gaps (e.g.\ $11.60\%$ versus $7.71\%$ at $(6,5,2)$, and $16.01\%$ versus $8.96\%$ at $(6,5,4)$), indicating that regularized cuts are weaker than normalized cuts in those cases when we alternate.

In contrast to the GEP and ARM instances, the alternating criterion provides little benefit on the CLSP instances. Solve times with the alternating criterion are often similar or higher than the default criterion: for example, Norm-D takes $293$s at $(3,3,8)$, but Norm-A takes $1016$s, and at $(6,4,8)$ Norm-D finishes in $3321$s with an $8.60\%$ gap while Norm-A hits the time limit with an $8.64\%$ gap. This behavior is explained by the Prop.\ column, which ranges from $0.36$ to $0.69$: roughly half of all iterations still require a ReLU cut because the Benders cut fails to separate the incumbent. The Alt~Gap column shows why: the LP--MIP subproblem gap ranges from $36\%$ to $69\%$ across all configurations. This is far larger than the gaps observed on most GEP ($\le 3\%$) and ARM ($\le 24\%$) instances. When the LP relaxation is a poor approximation of the MIP subproblem, the Benders cuts derived from LP duals are weak and may not cut off the incumbent as often, so the overhead of evaluating the MIP subproblem costs at every iteration yields no benefit, and the algorithm is forced to solve the expensive ReLU dual almost half the time, while producing weaker cuts overall. The large gaps on these instances, as well as the gaps observed on the ARM instances, can be partially attributed to bound stalling on a subset of scenario paths sampled during forward and backward passes of the SDDP algorithm. The compromise policy of \citealt{sen2016mitigating} and \citealt{xu2023compromise} addresses this by taking decisions based on the average of the approximate value functions learned via cutting planes across disjoint subsets of scenario paths, where each subset traverses a different set of scenario paths during the forward and backward passes. Averaging across different subsets in this way reduces variance and provides a tighter lower-bound estimate, offering a complementary approach for gap reduction on hard instances such as these.

\subsection{Portfolio Optimization}

\input{tables/paper/portfolio_test.tex}

Table~\ref{tab:portfolio_test} reports results for the portfolio optimization instances. The main observation is that the alternating criterion (A) yields a significant improvement across all configurations, reducing solve times from $721$--$3600$s under the default criterion (D) to $12$--$21$s, while simultaneously reducing the final gaps. This improvement is explained by the Prop.\ and Alt~Gap columns: Prop.\ $= 0.000$ implies that the Benders cut separated the incumbent in every iteration, and the ReLU dual was never invoked. The Alt~Gap is also $0.00\%$ across all configurations, confirming that the LP relaxation is exact for the portfolio subproblems.

Under the default criterion (D), Norm-D consistently achieves smaller gaps than Reg-D. For $N=25$ this advantage is moderate: $1.33\%$ vs.\ $2.49\%$ at $(N,T,R)=(25,6,5)$, $0.70\%$ vs.\ $3.97\%$ at $(25,6,10)$, and $1.76\%$ vs.\ $3.02\%$ at $(25,8,5)$. For $N=50$ the advantage becomes notable: Reg-D produces gaps of $32.10\%$, $55.32\%$, $45.45\%$, and $62.04\%$ at $(50,6,5)$, $(50,6,10)$, $(50,8,5)$, and $(50,8,10)$ respectively, versus $4.34\%$, $5.42\%$, $9.25\%$, and $9.71\%$ for Norm-D. These large gaps for Reg-D are similar to the trends seen on ARM, where regularized cuts are tight at the incumbent but weak globally, and the problem worsens as the state space grows larger with more assets. Norm-D is also faster in most configurations: average dual solution time (Dual t) is $1.97$s vs.\ $7.08$s at $(25,6,10)$, and $2.46$s vs.\ $5.16$s at $(25,6,5)$, consistent with the pattern seen on GEP and CLSP.

Under the alternating criterion, Norm-A and Reg-A are identical in every row with the same upper bound, gap, and iteration count, because neither ever invokes the ReLU dual (Prop.\ $= 0$).

\section{Conclusion}\label{sec:conclusion}

We presented an open-source software package, built on top of the SDDP.jl package (\citealt{dowson2021sddp}) in Julia, that implements ReLU cuts within a Benders-type decomposition framework for multistage stochastic integer programs, along with three cut-strengthening strategies: normalization, regularization, and a cut-generation criterion that alternates between cheaper Benders cuts and expensive tight cuts. To our knowledge, this is the first open-source implementation offering asymptotic convergence guarantees and cut-generation strategies for general MSIPs with mixed-integer state variables. Testing these strategies on four problem classes reveals several interesting trends. First, we observe that normalization is safer and more robust than regularization. Under the default cut-generating approach, normalization's final gaps are comparable to regularization's on GEP and CLSP but substantially smaller on ARM and portfolio, and its dual is also computationally cheaper to solve. Under the alternating criterion, both methods tend to produce cuts of similar quality, differing mainly in per-iteration cost rather than bound quality, with the notable exception of GEP, where normalization still yields consistently smaller gaps. The benefit of the alternating criterion depends on the LP--MIP subproblem gap. It yields an order-of-magnitude speedup on GEP, ARM, and portfolio, where this gap is small, but offers little benefit on CLSP, where the gap is large, and Benders cuts often fail to cut off the incumbent.

Based on these findings, we offer the following guidance for researchers and practitioners. When the subproblems are mixed-integer programs with general integer state variables, ReLU cuts are the appropriate choice, as they provide asymptotic convergence guarantees in this setting; and for binary state variables, Lagrangian cuts suffice at lower computational cost. Between the two cut-strengthening strategies, normalization tends to be the safer default: it consistently produces stronger bounds across all four problem classes, whereas regularization may fail to improve on trivial bounds when the state space is large, as seen on the ARM and portfolio instances. Both strategies are sensitive to the choice of core point, so we recommend using the tune-then-test protocol described in Section~\ref{sec:comp} to select the core point family and its scaling parameters on a held-out instance before evaluating on test data. Finally, the alternating criterion is highly effective when the LP--MIP subproblem gap is small, as cheaper Benders cuts then separate the incumbent in most iterations and the expensive ReLU dual is invoked rarely. When the LP--MIP gap is large, as in the CLSP instances, Benders cuts fail frequently, and the alternating criterion provides little benefit; in this regime, the default (non-alternating) configuration is preferable.

\appendix

\section{Problem Classes}\label{sec:appendix}

In this section, we describe the problem classes in greater detail, including their uncertain parameters, state and local variables, constraints, and objective functions.

\subsection{Generation Expansion Planning}

We first discuss the generation expansion planning (GEP) problem (\citealt{zou2019stochastic, jin2011modeling}). There are $n$ types of generators that are available. Let $g_t$ and $x_t$ denote the decision vectors indicating the number of generators of each type built at stage $t$, and the cumulative number of generators of each type built until stage $t$, respectively. Each time period $t$ is further divided into sub-periods. Following \citealt{jin2011modeling}, we consider three sub-periods indicating low, medium, and high demand. Let $y_{ts}$ be the local variable vector denoting the electricity produced by each generator in time period $t$ and sub-period $s$. The demand in stage $t$ and sub-period $s$ is denoted by $d_{ts}$ while the unmet demand is given by $u_{ts}$. Now, the optimization problem at stage $t$ is given by:
\begin{subequations}
\begin{align}
    Q_t(x_{t-1}, w_t^j) := & \min_{g_t, y_{ts}} \hspace{0.2cm} a_t^{\top} g_t + \sum_{s} (q_s(w_t^j))^{\top} y_{ts}+ \rho\sum_{s} u_{ts} + \sum_{k=1}^{N_{t+1}}q_k Q_{t+1}(x_t, w_{t+1}^k)  \\ 
    &x_t \leq G \label{eq:gep_bound}\\ 
    &x_t = x_{t-1} + g_t \label{eq:gep_connect}\\
    &y_{ts} \leq A_t x_t, &s \in S \label{eq:gep_elec}\\
    &\mathbf{1}^{\top} y_{ts} + u_{ts} = d_{ts}(w_t^j), &s \in S \label{eq:gep_dem}\\
    &x_t, g_t \in \mathbb{Z}_{+}^{n}, y_t \in \mathbb{R}^{n}_{+}. \nonumber
\end{align}
\end{subequations}

The uncertainty appears in the objective function through the coefficient vector $q_s(w_t^j)$. It also appears in constraints because demand $d_{ts}$ is random. The uncertainty realizations are stagewise independent, and the subproblem depends only on the stage, not on the specific node. Therefore, we replace the node-based notation $Q_n$ in \eqref{prob:sub} with the stage-based notation $Q_t$. The constraint \eqref{eq:gep_bound} enforces that only a pre-determined number of generators $G$ can be used across all time periods. In constraint \eqref{eq:gep_connect}, the decision variables in the current stage connect with those in previous stages. The constraint \eqref{eq:gep_elec} limits the generation capacity based on the available number of generators. The matrix $A_t$ contains the maximum rating and maximum capacity information of generators. Finally, the constraint \eqref{eq:gep_dem} enforces that the demand is met. The objective function evaluates the cost of adding generators and generating electricity while penalizing unmet demand. The first summand captures the cost of building the generators via the coefficients $a_t$, the second summand is the cost of generating electricity, and the third summand is the penalty cost for failing to meet demand. The final summand is the expected future cost arising from stage $t+1$ onwards.

\subsection{Airline Revenue Management}

We consider a multistage stochastic network revenue management problem over $T$ booking intervals, following the model in \citealt{moller2008airline}. The objective is to maximize expected ticket revenue from seat sales on a flight network. Each component of the $m$-dimensional decision vector corresponds to an origin--destination itinerary and fare-class pair. We need to decide how to allocate limited seats across different itineraries and fare classes, given uncertain passenger demand.

Let $t \in \{1,\ldots,T\}$ index booking intervals, $\mathcal{I}$ the set of itineraries, and $\mathcal{J}$ the set of fare classes. Each pair $(i,j) \in \mathcal{I} \times \mathcal{J}$ maps to one component of the $m$-vector. Legs and cabins index the capacity constraints. The random parameter at stage $t$ is the demand vector $d_t \in \mathbb{Z}_+^m$, whose components bound the number of booking requests for each itinerary--fare-class pair. Scenarios are generated using the same data-generation procedure in \citealt{zou2019stochastic}.

The integer state variables carried across stages are $B_t \in \mathbb{Z}_+^m$, cumulative fulfilled bookings through interval $t$, and $C_t \in \mathbb{Z}_+^m$, cumulative cancellations through interval $t$, with $B_0 = C_0 = 0$. At each stage $t$, the local integer decision variables are
$b_t \in \mathbb{Z}_+^m$, new bookings in interval $t$, and
$c_t \in \mathbb{Z}_+^m$, new cancellations in interval $t$.

The optimization problem at stage $t$ is given by:
\begin{subequations}\label{eq:arm}
\begin{align}
    Q_t(B_{t-1}, C_{t-1}, \omega_t) := & \max_{b_t, c_t, B_t, C_t} \hspace{0.2cm} (f^b)^\top b_t - (f^c)^\top c_t + \sum_{k=1}^{N_{t+1}} q_k\, Q_{t+1}(B_t, C_t, \omega_{t+1}^k) \notag \\
    & B_t = B_{t-1} + b_t, \quad C_t = C_{t-1} + c_t, \label{eq:arm-flow}\\
    & C_t = \lfloor \Phi_t B_t + \tfrac{1}{2} \rfloor, \label{eq:arm-crate}\\
    & b_t \le d_t(\omega_t), \label{eq:arm-demand}\\
    & A(B_t - C_t) \le R, \label{eq:arm-capacity}\\
    & B_t,\, C_t,\, b_t,\, c_t \in \mathbb{Z}_+^m. \nonumber
\end{align}
\end{subequations}
Here, $f^b$ and $f^c$ denote the booking price and cancellation refund vectors, respectively. Constraint \eqref{eq:arm-flow} tracks cumulative bookings and cancellations across stages. The cancellation rate constraint \eqref{eq:arm-crate} links cancellations to bookings via the diagonal matrix $\Phi_t$ of fare-class cancellation rates. Constraint \eqref{eq:arm-demand} caps new bookings by the realized demand $d_t(\omega_t)$, and \eqref{eq:arm-capacity} enforces seat limits via the $0$--$1$ incidence matrix $A$, which indicates whether an itinerary--fare-class pair consumes capacity on a given leg and cabin, and the seat-limit vector $R$.

\subsection{Capacitated Lot-Sizing}

We consider a multistage stochastic capacitated lot-sizing problem (CLSP) on a single shared production resource, following the setup-time model of \citealt{trigeiro1989capacitated}. A firm plans production for $P$ products over $T$ stages. At each stage $t$, the uncertainty is a random vector $\omega_t = (\xi_{t,1}, \ldots, \xi_{t,P})$ of multiplicative demand factors. The realized demand for product $i \in \mathcal{P} := \{1,\ldots,P\}$ is $d_{t,i}(\omega_t) = \xi_{t,i}\,\bar{d}_{t,i}$, where $\bar{d}_{t,i}$ is a deterministic mean-demand parameter. Each component $\xi_{t,i}$ is drawn from a product-specific lognormal distribution; stage $1$ is deterministic, and at stages $t = 2,\ldots,T$ we use $R$ equally likely realizations. We follow the same data generation procedure as in \citealt{fullnernew}. The state variable $I_{t,i}$ denotes end-of-stage inventory for product $i$ after stage $t$, with $I_{0,i} = 0$ for all $i$. At each stage, the local decisions are $x_{t,i} \ge 0$ (production quantity), $l_{t,i} \ge 0$ (lost sales), and $y_{t,i} \in \{0,1\}$ (setup indicator) for each product $i \in \mathcal{P}$. The fixed parameters are $s_i$ (setup time), $a_i$ (production time per unit), $f_i$ (setup cost), $h_i$ (inventory holding cost), $p_i$ (lost-sales penalty), $C$ (shared capacity in time units), and $\bar{I}$ (inventory upper bound).

The optimization problem at stage $t$ is given by:
\begin{subequations}\label{clsp:constr}
\begin{align}
    Q_t(I_{t-1}, \omega_t) := & \min_{x_{t}, l_{t}, y_{t}, I_t} \hspace{0.2cm} \sum_{i \in \mathcal{P}} \bigl(f_i y_{t,i} + h_i I_{t,i} + p_i l_{t,i}\bigr) + \sum_{k=1}^{N_{t+1}} q_k\, Q_{t+1}(I_t, \omega_{t+1}^k) \notag \\
    & x_{t,i} \le C\, y_{t,i}, && \forall i \in \mathcal{P}, \label{clsp:setup} \\
    & \sum_{i \in \mathcal{P}} \bigl(a_i x_{t,i} + s_i y_{t,i}\bigr) \le C, \label{clsp:cap} \\
    & I_{t,i} - l_{t,i} = I_{t-1,i} + x_{t,i} - d_{t,i}(\omega_t), && \forall i \in \mathcal{P}, \label{clsp:balance} \\
    & 0 \le I_{t,i} \le \bar{I}, \quad x_{t,i},\, l_{t,i} \ge 0, && \forall i \in \mathcal{P}, \nonumber \\
    & y_{t,i} \in \{0,1\}, && \forall i \in \mathcal{P}. \nonumber
\end{align}
\end{subequations}
Constraint~\eqref{clsp:setup} permits production for product $i$ only if a setup is performed (big-$M$ linking with $M = C$). Constraint~\eqref{clsp:cap} is the shared capacity constraint, i.e., setup times and production times compete for the same time budget $C$. Constraint~\eqref{clsp:balance} is the inventory balance with lost sales, where ending inventory net of lost sales equals beginning inventory plus production minus realized demand $d_{t,i}(\omega_t)$. The objective minimizes the total setup, holding, and lost-sales costs at stage $t$, plus the expected future costs from stage $t+1$ onwards. The instances are taken directly from \citealt{fullnernew} and do not include production cost in the objective. 

\subsection{Portfolio Optimization}

We consider a multistage portfolio optimization problem from \citealt{zou2019stochastic}. An investor holds cash (asset $c := 1$) and $N-1$ risky stocks (indexed by $I := \{2,\ldots,N\}$) over $T$ stages. At each stage, the investor may buy or sell stocks, subject to proportional transaction costs $\alpha_b$ and $\alpha_s$ for purchases and sales, respectively, and under a cardinality constraint that limits the number of stocks held simultaneously to at most $K$. The uncertainty at stage $t$ is the return vector $\omega_t \in \mathbb{R}^N$, where $\omega_{t,i}$ is the gross return factor for asset $i$. Returns are drawn from historical bi-weekly S\&P data; stage $1$ is deterministic, and at stages $t = 2,\ldots,T$ we sample $R$ return vectors independently with replacement, each with probability $1/R$. The state variable $x_{t,i}$ denotes the end-of-stage dollar value held in asset $i$ after trades at stage $t$, with $x_{0,c} = \bar{x}_0$ and $x_{0,i} = 0$ for $i \in I$. The local decisions at each stage are $b_{t,i} \ge 0$ (dollar amount bought in stock $i$), $s_{t,i} \ge 0$ (dollar amount sold from stock $i$), and $z_{t,i} \in \{0,1\}$ (indicator that stock $i$ is held), for each $i \in I$. The parameters $M$, $v$, and $u$ are big-$M$, per-asset value, and per-trade upper bounds, respectively.

The optimization problem at stage $t$ is given by:
\begin{subequations}\label{por:constr}
\begin{align}
    Q_t(x_{t-1}, \omega_t) := & \max_{b_t, s_t, z_t, x_t} \hspace{0.2cm} \sum_{k=1}^{N_{t+1}} q_k\, Q_{t+1}(x_t, \omega_{t+1}^k) \notag \\
    & x_{t,i} = \omega_{t,i}\, x_{t-1,i} + b_{t,i} - s_{t,i}, && \forall i \in I, \label{por:flow} \\
    & x_{t,c} = \omega_{t,c}\, x_{t-1,c} - (1+\alpha_b) \textstyle\sum_{i \in I} b_{t,i} + (1-\alpha_s) \textstyle\sum_{i \in I} s_{t,i}, \label{por:cash} \\
    & x_{t,i} \le M\, z_{t,i}, \quad s_{t,i} \le \omega_{t,i}\, x_{t-1,i}, && \forall i \in I, \label{por:hold} \\
    & \textstyle\sum_{i \in I} z_{t,i} \le K, \label{por:card} \\
    & 0 \le b_{t,i},\, s_{t,i} \le u, \quad 0 \le x_{t,i} \le v, && \forall i \in I, \notag \\
    & 0 \le x_{t,c} \le v, \quad z_{t,i} \in \{0,1\}, && \forall i \in I. \notag
\end{align}
\end{subequations}
The objective at each intermediate stage $t < T$ has no immediate reward; the terminal value function $Q_T(x_{T-1}, \omega_T)$ returns $\sum_{i=1}^{N} x_{T,i}$, the total end-of-horizon wealth. Constraint~\eqref{por:flow} is the stock transaction balance and \eqref{por:cash} is the self-financing cash balance, where purchases deplete cash at rate $(1+\alpha_b)$ and sales replenish it at rate $(1-\alpha_s)$. Constraints~\eqref{por:hold} enforce big-$M$ linking of holdings to the binary indicators and cap sales at the post-return value of beginning holdings. Constraint~\eqref{por:card} limits the number of stocks held to at most $K$.

\bibliographystyle{plainnat}
\bibliography{ref.bib}

\end{document}

%% file: tables/paper/gep_test.tex
{\scriptsize
\setlength{\tabcolsep}{3pt}
\begin{longtable}[t]{cccccccccccc}
\caption{Out-of-sample test performance on multistage GEP instances.}\label{tab:gep_test}\\
\toprule
$T$ & $R$ & Type & Method & LB & Time (s) & Iters & Dual-Iters & Prop. & Dual t (s) & Alt Gap (\%) & Gap (\%) \\
\midrule
\endfirsthead
\multicolumn{12}{l}{\small\itshape (continued from previous page)}\\
\toprule
$T$ & $R$ & Type & Method & LB & Time (s) & Iters & Dual-Iters & Prop. & Dual t (s) & Alt Gap (\%) & Gap (\%) \\
\midrule
\endhead
\midrule
\multicolumn{12}{r}{\small\itshape (continued on next page)}\\
\endfoot
\bottomrule
\endlastfoot
10 & 3 & D & Norm & 24435.21 & 214 & 18 & 29 &  & 0.46 &  & 1.11 \\
 &  &  & Reg & 24469.62 & 187 & 16 & 41 &  & 0.50 &  & 0.95 \\
 &  & A & Norm & 24204.83 & 26 & 23 & 30 & 0.016 & 0.31 & 1.61 & 2.57 \\
 &  &  & Reg & 24150.71 & 20 & 22 & 14 & 0.015 & 0.14 & 1.49 & 4.50 \\
\addlinespace
10 & 6 & D & Norm & 26176.84 & 504 & 19 & 28 &  & 0.54 &  & 1.15 \\
 &  &  & Reg & 26219.96 & 665 & 16 & 43 &  & 0.87 &  & 0.88 \\
 &  & A & Norm & 25902.26 & 38 & 21 & 33 & 0.025 & 0.44 & 2.52 & 2.46 \\
 &  &  & Reg & 25895.70 & 20 & 21 & 8 & 0.025 & 0.07 & 2.47 & 2.58 \\
\addlinespace
10 & 10 & D & Norm & 24574.06 & 776 & 17 & 30 &  & 0.57 &  & 1.06 \\
 &  &  & Reg & 24577.24 & 1536 & 18 & 44 &  & 1.04 &  & 1.02 \\
 &  & A & Norm & 24351.54 & 42 & 25 & 30 & 0.024 & 0.37 & 2.41 & 1.88 \\
 &  &  & Reg & 24188.23 & 28 & 23 & 15 & 0.030 & 0.12 & 3.00 & 3.19 \\
\addlinespace
10 & 12 & D & Norm & 24982.07 & 1347 & 18 & 30 &  & 0.75 &  & 1.12 \\
 &  &  & Reg & 25011.11 & 1723 & 18 & 44 &  & 1.04 &  & 1.02 \\
 &  & A & Norm & 24720.79 & 55 & 22 & 34 & 0.019 & 0.49 & 1.91 & 2.44 \\
 &  &  & Reg & 24637.17 & 28 & 21 & 7 & 0.019 & 0.07 & 1.94 & 5.68 \\
\addlinespace
10 & 20 & D & Norm & 24690.80 & 2334 & 18 & 29 &  & 0.76 &  & 1.24 \\
 &  &  & Reg & 24716.84 & 3589 & 18 & 45 &  & 1.26 &  & 1.05 \\
 &  & A & Norm & 24470.50 & 144 & 21 & 29 & 0.033 & 0.95 & 3.28 & 2.07 \\
 &  &  & Reg & 24382.94 & 30 & 21 & 7 & 0.032 & 0.06 & 3.16 & 3.62 \\
\addlinespace
10 & 50 & D & Norm & 24265.09 & TL & 13 & 25 &  & 0.70 &  & 2.35 \\
 &  &  & Reg & 24333.99 & TL & 8 & 43 &  & 1.07 &  & 1.93 \\
 &  & A & Norm & 24168.05 & 212 & 21 & 29 & 0.029 & 0.78 & 2.91 & 3.72 \\
 &  &  & Reg & 24222.02 & 108 & 22 & 12 & 0.029 & 0.23 & 2.86 & 3.14 \\
\addlinespace
15 & 3 & D & Norm & 76519.07 & 196 & 15 & 25 &  & 0.33 &  & 3.31 \\
 &  &  & Reg & 76573.40 & 595 & 18 & 43 &  & 0.94 &  & 3.12 \\
 &  & A & Norm & 75920.42 & 22 & 18 & 23 & 0.003 & 0.27 & 0.29 & 4.88 \\
 &  &  & Reg & 75826.33 & 19 & 17 & 14 & 0.003 & 0.18 & 0.25 & 4.58 \\
\addlinespace
20 & 3 & D & Norm & 1197665.18 & 129 & 15 & 26 &  & 0.15 &  & 0.85 \\
 &  &  & Reg & 1179358.19 & TL & 40 & 36 &  & 1.78 &  & 2.80 \\
 &  & A & Norm & 1196634.62 & 19 & 11 & 51 & 0.000 & 0.42 & 0.00 & 1.01 \\
 &  &  & Reg & 1196562.08 & 19 & 11 & 12 & 0.000 & 0.61 & 0.00 & 1.08 \\
\end{longtable}
}

%% file: tables/paper/arm_test.tex
{\scriptsize
\setlength{\tabcolsep}{3pt}
\begin{longtable}[t]{cccccccccccc}
\caption{Out-of-sample test performance on multistage airline revenue management (ARM) instances. UB is an upper bound on the (maximization) objective.}\label{tab:arm_test}\\
\toprule
$T$ & $R$ & Type & Method & UB & Time (s) & Iters & Dual-Iters & Prop. & Dual t (s) & Alt Gap (\%) & Gap (\%) \\
\midrule
\endfirsthead
\multicolumn{12}{l}{\small\itshape (continued from previous page)}\\
\toprule
$T$ & $R$ & Type & Method & UB & Time (s) & Iters & Dual-Iters & Prop. & Dual t (s) & Alt Gap (\%) & Gap (\%) \\
\midrule
\endhead
\midrule
\multicolumn{12}{r}{\small\itshape (continued on next page)}\\
\endfoot
\bottomrule
\endlastfoot
4 & 2 & D & Norm & 228776.77 & TL & 10 & 113 &  & 63.62 &  & 14.16 \\
 &  &  & Reg & 1059840.00 & 465 & 11 & 43 &  & 10.12 &  & 82.98 \\
 &  & A & Norm & 214085.08 & 16 & 19 & 25 & 0.018 & 1.12 & 22.37 & 3.41 \\
 &  &  & Reg & 213859.95 & 134 & 20 & 294 & 0.065 & 14.14 & 7.79 & 3.25 \\
\addlinespace
4 & 4 & D & Norm & 245958.02 & TL & 8 & 97 &  & 39.63 &  & 20.95 \\
 &  &  & Reg & 1059840.00 & 877 & 11 & 42 &  & 9.70 &  & 82.66 \\
 &  & A & Norm & 214014.77 & 21 & 19 & 18 & 0.014 & 1.23 & 0.40 & 3.77 \\
 &  &  & Reg & 214014.77 & 243 & 19 & 300 & 0.055 & 19.99 & 0.42 & 3.77 \\
\addlinespace
4 & 6 & D & Norm & 249358.58 & TL & 6 & 58 &  & 36.96 &  & 22.20 \\
 &  &  & Reg & 1059840.00 & 1910 & 11 & 40 &  & 14.21 &  & 82.34 \\
 &  & A & Norm & 214118.99 & 25 & 19 & 17 & 0.015 & 1.60 & 0.43 & 3.74 \\
 &  &  & Reg & 213930.73 & 371 & 19 & 300 & 0.034 & 28.97 & 0.42 & 3.63 \\
\addlinespace
5 & 2 & D & Norm & 235087.31 & TL & 14 & 55 &  & 38.78 &  & 16.25 \\
 &  &  & Reg & 1059840.00 & 269 & 11 & 24 &  & 5.61 &  & 82.02 \\
 &  & A & Norm & 214547.90 & 21 & 18 & 31 & 0.019 & 2.48 & 14.62 & 4.86 \\
 &  &  & Reg & 214535.63 & 82 & 18 & 163 & 0.031 & 11.10 & 14.08 & 4.86 \\
\addlinespace
5 & 4 & D & Norm & 242154.22 & TL & 9 & 50 &  & 30.54 &  & 20.30 \\
 &  &  & Reg & 1059840.00 & 370 & 11 & 21 &  & 4.10 &  & 82.21 \\
 &  & A & Norm & 214811.51 & 42 & 20 & 38 & 0.012 & 5.15 & 7.40 & 4.66 \\
 &  &  & Reg & 214789.64 & 119 & 20 & 300 & 0.012 & 27.07 & 7.39 & 4.69 \\
\addlinespace
5 & 6 & D & Norm & 252042.17 & TL & 7 & 48 &  & 26.93 &  & 22.07 \\
 &  &  & Reg & 1059840.00 & 866 & 11 & 26 &  & 6.99 &  & 82.13 \\
 &  & A & Norm & 216673.31 & 30 & 19 & 31 & 0.004 & 6.72 & 19.42 & 4.27 \\
 &  &  & Reg & 216670.62 & 79 & 19 & 227 & 0.005 & 23.54 & 19.41 & 4.27 \\
\addlinespace
6 & 2 & D & Norm & 257241.31 & TL & 8 & 56 &  & 55.27 &  & 24.33 \\
 &  &  & Reg & 1059840.00 & 212 & 11 & 22 &  & 4.31 &  & 82.49 \\
 &  & A & Norm & 216196.19 & 34 & 20 & 24 & 0.017 & 3.24 & 7.1e8 & 5.55 \\
 &  &  & Reg & 216130.99 & 70 & 20 & 239 & 0.018 & 17.21 & 7.1e8 & 5.56 \\
\addlinespace
6 & 4 & D & Norm & 252458.02 & TL & 7 & 50 &  & 28.24 &  & 23.66 \\
 &  &  & Reg & 1059840.00 & 364 & 11 & 21 &  & 4.16 &  & 82.41 \\
 &  & A & Norm & 213949.11 & 38 & 21 & 30 & 0.013 & 2.96 & 18.10 & 4.55 \\
 &  &  & Reg & 214005.76 & 119 & 20 & 234 & 0.015 & 19.34 & 15.25 & 4.58 \\
\addlinespace
6 & 6 & D & Norm & 258714.52 & TL & 5 & 50 &  & 30.99 &  & 25.69 \\
 &  &  & Reg & 1059840.00 & 680 & 11 & 22 &  & 5.14 &  & 82.58 \\
 &  & A & Norm & 214639.46 & 33 & 21 & 20 & 0.007 & 2.65 & 0.89 & 4.65 \\
 &  &  & Reg & 214592.21 & 202 & 21 & 300 & 0.010 & 27.95 & 0.66 & 4.68 \\
\end{longtable}
}

%% file: tables/paper/clsp_test.tex
{\scriptsize
\setlength{\tabcolsep}{3pt}
\begin{longtable}[t]{ccccccccccccc}
\caption{Out-of-sample test performance on multistage CLSP instances. The column $P$ in table denotes the number of products considered in the problem.}\label{tab:clsp_test}\\
\toprule
$P$ & $T$ & $R$ & Type & Method & LB & Time (s) & Iters & Dual-Iters & Prop. & Dual t (s) & Alt Gap (\%) & Gap (\%) \\
\midrule
\endfirsthead
\multicolumn{13}{l}{\small\itshape (continued from previous page)}\\
\toprule
$P$ & $T$ & $R$ & Type & Method & LB & Time (s) & Iters & Dual-Iters & Prop. & Dual t (s) & Alt Gap (\%) & Gap (\%) \\
\midrule
\endhead
\midrule
\multicolumn{13}{r}{\small\itshape (continued on next page)}\\
\endfoot
\bottomrule
\endlastfoot
3 & 3 & 8 & D & Norm & 1046.66 & 293 & 25 & 24 &  & 0.79 &  & 7.84 \\
 &  &  &  & Reg & 1045.65 & 304 & 28 & 27 &  & 0.78 &  & 8.84 \\
 &  &  & A & Norm & 1048.25 & 1016 & 36 & 23 & 0.523 & 2.57 & 52.33 & 7.63 \\
 &  &  &  & Reg & 1048.83 & 468 & 31 & 32 & 0.539 & 1.32 & 53.90 & 9.04 \\
\addlinespace
3 & 4 & 4 & D & Norm & 1732.03 & 189 & 29 & 22 &  & 0.55 &  & 6.09 \\
 &  &  &  & Reg & 1727.44 & 223 & 29 & 24 &  & 0.65 &  & 6.75 \\
 &  &  & A & Norm & 1738.50 & 263 & 37 & 20 & 0.432 & 0.79 & 43.17 & 5.71 \\
 &  &  &  & Reg & 1737.43 & 269 & 33 & 31 & 0.424 & 1.13 & 42.43 & 5.60 \\
\addlinespace
3 & 4 & 8 & D & Norm & 1710.70 & 551 & 29 & 22 &  & 0.84 &  & 4.64 \\
 &  &  &  & Reg & 1704.33 & 368 & 23 & 28 &  & 0.72 &  & 5.15 \\
 &  &  & A & Norm & 1713.68 & 640 & 37 & 23 & 0.423 & 0.99 & 42.33 & 4.40 \\
 &  &  &  & Reg & 1708.96 & 818 & 39 & 26 & 0.424 & 1.27 & 42.40 & 7.14 \\
\addlinespace
3 & 5 & 2 & D & Norm & 1809.22 & 134 & 29 & 21 &  & 0.54 &  & 7.42 \\
 &  &  &  & Reg & 1815.41 & 173 & 28 & 32 &  & 0.80 &  & 7.87 \\
 &  &  & A & Norm & 1810.92 & 138 & 36 & 24 & 0.387 & 0.57 & 38.69 & 7.43 \\
 &  &  &  & Reg & 1813.63 & 170 & 37 & 33 & 0.366 & 0.81 & 36.60 & 7.82 \\
\addlinespace
3 & 5 & 4 & D & Norm & 2005.65 & 344 & 29 & 22 &  & 0.73 &  & 4.70 \\
 &  &  &  & Reg & 2004.42 & 713 & 31 & 32 &  & 1.39 &  & 3.95 \\
 &  &  & A & Norm & 2004.61 & 658 & 47 & 24 & 0.357 & 1.41 & 35.66 & 3.87 \\
 &  &  &  & Reg & 2009.93 & 1124 & 46 & 34 & 0.347 & 2.24 & 34.74 & 3.92 \\
\addlinespace
6 & 3 & 8 & D & Norm & 1865.31 & 2982 & 35 & 47 &  & 5.27 &  & 7.94 \\
 &  &  &  & Reg & 1825.40 & 3326 & 32 & 48 &  & 6.47 &  & 9.68 \\
 &  &  & A & Norm & 1873.79 & TL & 42 & 44 & 0.640 & 6.96 & 63.95 & 6.89 \\
 &  &  &  & Reg & 1793.09 & TL & 41 & 46 & 0.691 & 8.33 & 69.05 & 12.74 \\
\addlinespace
6 & 4 & 4 & D & Norm & 2710.43 & 2012 & 38 & 46 &  & 4.31 &  & 6.11 \\
 &  &  &  & Reg & 2704.54 & 2828 & 36 & 58 &  & 6.59 &  & 7.49 \\
 &  &  & A & Norm & 2720.80 & 2612 & 51 & 47 & 0.576 & 5.75 & 57.56 & 7.13 \\
 &  &  &  & Reg & 2655.76 & 3281 & 48 & 50 & 0.602 & 7.39 & 60.23 & 8.83 \\
\addlinespace
6 & 4 & 8 & D & Norm & 2767.28 & 3321 & 30 & 46 &  & 4.51 &  & 8.60 \\
 &  &  &  & Reg & 2766.39 & TL & 27 & 60 &  & 6.01 &  & 7.81 \\
 &  &  & A & Norm & 2733.21 & TL & 41 & 47 & 0.594 & 5.24 & 59.35 & 8.64 \\
 &  &  &  & Reg & 2600.98 & TL & 43 & 42 & 0.634 & 5.47 & 63.35 & 14.88 \\
\addlinespace
6 & 5 & 2 & D & Norm & 3393.80 & 725 & 32 & 45 &  & 2.85 &  & 8.76 \\
 &  &  &  & Reg & 3374.58 & 1330 & 37 & 57 &  & 4.86 &  & 12.06 \\
 &  &  & A & Norm & 3416.61 & 1060 & 48 & 45 & 0.478 & 3.74 & 47.77 & 7.71 \\
 &  &  &  & Reg & 3357.22 & 2113 & 57 & 52 & 0.502 & 6.48 & 50.25 & 11.60 \\
\addlinespace
6 & 5 & 4 & D & Norm & 3354.76 & 2281 & 33 & 46 &  & 4.46 &  & 6.60 \\
 &  &  &  & Reg & 3298.79 & TL & 36 & 55 &  & 6.81 &  & 10.14 \\
 &  &  & A & Norm & 3300.24 & 3544 & 49 & 47 & 0.501 & 6.44 & 50.13 & 8.96 \\
 &  &  &  & Reg & 3119.80 & TL & 44 & 58 & 0.521 & 9.94 & 52.14 & 16.01 \\
\end{longtable}
}

%% file: tables/paper/portfolio_test.tex
{\scriptsize
\setlength{\tabcolsep}{3pt}
\begin{longtable}[t]{ccccccccccccc}
\caption{Out-of-sample test performance on multistage portfolio optimization instances. UB is an upper bound on the (maximization) objective; Gap (\%) is $|\text{UB}-\text{LB}|/|\text{UB}|$, averaged over three test instances.}\label{tab:portfolio_test}\\
\toprule
$N$ & $T$ & $R$ & Type & Method & UB & Time (s) & Iters & Dual-Iters & Prop. & Dual t (s) & Alt Gap (\%) & Gap (\%) \\
\midrule
\endfirsthead
\multicolumn{13}{l}{\small\itshape (continued from previous page)}\\
\toprule
$N$ & $T$ & $R$ & Type & Method & UB & Time (s) & Iters & Dual-Iters & Prop. & Dual t (s) & Alt Gap (\%) & Gap (\%) \\
\midrule
\endhead
\midrule
\multicolumn{13}{r}{\small\itshape (continued on next page)}\\
\endfoot
\bottomrule
\endlastfoot
25 & 6 & 5 & D & Norm & 115.39 & 721 & 14 & 127 &  & 2.46 &  & 1.33 \\
 &  &  &  & Reg & 116.27 & 1303 & 15 & 163 &  & 5.16 &  & 2.49 \\
 &  &  & A & Norm & 114.53 & 13 & 15 &  & 0.000 &  & 0.00 & 0.36 \\
 &  &  &  & Reg & 114.53 & 13 & 15 &  & 0.000 &  & 0.00 & 0.36 \\
\addlinespace
25 & 6 & 10 & D & Norm & 114.38 & 799 & 15 & 114 &  & 1.97 &  & 0.70 \\
 &  &  &  & Reg & 116.16 & 3398 & 13 & 163 &  & 7.08 &  & 3.97 \\
 &  &  & A & Norm & 114.17 & 15 & 12 &  & 0.000 &  & 0.00 & 0.53 \\
 &  &  &  & Reg & 114.17 & 13 & 12 &  & 0.000 &  & 0.00 & 0.53 \\
\addlinespace
25 & 8 & 5 & D & Norm & 125.18 & 1446 & 19 & 123 &  & 2.81 &  & 1.76 \\
 &  &  &  & Reg & 125.83 & 863 & 14 & 146 &  & 2.95 &  & 3.02 \\
 &  &  & A & Norm & 124.00 & 14 & 16 &  & 0.000 &  & 0.00 & 0.51 \\
 &  &  &  & Reg & 124.00 & 12 & 16 &  & 0.000 &  & 0.00 & 0.51 \\
\addlinespace
25 & 8 & 10 & D & Norm & 116.15 & TL & 14 & 121 &  & 4.38 &  & 2.80 \\
 &  &  &  & Reg & 121.16 & 3394 & 12 & 168 &  & 5.23 &  & 9.18 \\
 &  &  & A & Norm & 114.00 & 14 & 19 &  & 0.000 &  & 0.00 & 1.33 \\
 &  &  &  & Reg & 114.00 & 21 & 19 &  & 0.000 &  & 0.00 & 1.33 \\
\addlinespace
50 & 6 & 5 & D & Norm & 119.60 & 1710 & 15 & 181 &  & 4.56 &  & 4.34 \\
 &  &  &  & Reg & 153.53 & 945 & 12 & 187 &  & 3.99 &  & 32.10 \\
 &  &  & A & Norm & 116.47 & 16 & 16 &  & 0.000 &  & 0.00 & 1.34 \\
 &  &  &  & Reg & 116.47 & 16 & 16 &  & 0.000 &  & 0.00 & 1.34 \\
\addlinespace
50 & 6 & 10 & D & Norm & 118.18 & 2735 & 13 & 179 &  & 4.87 &  & 5.42 \\
 &  &  &  & Reg & 228.33 & 1783 & 11 & 189 &  & 4.34 &  & 55.32 \\
 &  &  & A & Norm & 114.63 & 17 & 16 &  & 0.000 &  & 0.00 & 0.53 \\
 &  &  &  & Reg & 114.63 & 16 & 16 &  & 0.000 &  & 0.00 & 0.53 \\
\addlinespace
50 & 8 & 5 & D & Norm & 132.19 & 2472 & 12 & 174 &  & 7.04 &  & 9.25 \\
 &  &  &  & Reg & 187.25 & 2275 & 13 & 188 &  & 5.38 &  & 45.45 \\
 &  &  & A & Norm & 124.33 & 16 & 17 &  & 0.000 &  & 0.00 & 2.72 \\
 &  &  &  & Reg & 124.33 & 14 & 17 &  & 0.000 &  & 0.00 & 2.72 \\
\addlinespace
50 & 8 & 10 & D & Norm & 128.96 & TL & 11 & 171 &  & 5.40 &  & 9.71 \\
 &  &  &  & Reg & 275.15 & 3402 & 12 & 182 &  & 5.12 &  & 62.04 \\
 &  &  & A & Norm & 120.70 & 18 & 16 &  & 0.000 &  & 0.00 & 0.62 \\
 &  &  &  & Reg & 120.70 & 17 & 16 &  & 0.000 &  & 0.00 & 0.62 \\
\end{longtable}
}